\documentclass[12pt]{article}
\usepackage{amsfonts}
\usepackage{amssymb}
\usepackage{graphicx}
\usepackage{amsmath}

\begin{document}

\begin{center}
{\Large On The Problem Of Entire Functions That Share a Small Function With
Their Difference Operators}

\quad

\textbf{Abdallah} \textbf{EL FARISSI}$^{1}$, \textbf{Zinel\^{a}abidine}
\textbf{LATREUCH}$^{2}$, \textbf{\ Benharrat} \textbf{BELA\"{I}DI}$^{2}$
\textbf{and Asim ASIRI}$^{3}$

\quad

$^{1}$\textbf{Department of Mathematics and Informatics, }

\textbf{Faculty of Exact Sciences,}

\textbf{University of Bechar-(Algeria)}

\textbf{elfarissi.abdallah@yahoo.fr}

$^{2}$\textbf{Department of Mathematics }

\textbf{Laboratory of Pure and Applied Mathematics }

\textbf{University of Mostaganem (UMAB) }

\textbf{B. P. 227 Mostaganem-(Algeria)}

\textbf{z.latreuch@gmail.com}

\textbf{belaidibenharrat@yahoo.fr}

$^{3}$\textbf{Department of Mathematics, }

\textbf{Faculty of Science, King Abdulaziz University, }

\textbf{P.O.Box 80203, Jeddah 21589, Saudi Arabia}

\textbf{amkasiri@kau.edu.sa}

\quad
\end{center}

\noindent \textbf{Abstract. }In this paper, we study uniqueness problems for
an entire function that shares small functions of finite order with their
difference operators. In particular, we give a generalization of results in $%
[2,3,13]$.

\quad



\section{Introduction and Main Results}

\noindent In this paper, by meromorphic functions we will always mean
meromorphic functions in the complex plane. In what follows, we assume that
the reader is familiar with the fundamental results and the standard
notations of the Nevanlinna's value distribution theory of meromorphic
functions $(\left[ 9\right] ,$ $\left[ 11\right] ,$ $\left[ 17\right] )$. In
addition, we will use $\rho \left( f\right) $ to denote the order of growth
of $f$ and $\lambda \left( f\right) $ to denote the exponent of convergence
of zeros of $f$, we say that a meromorphic function $\varphi \left( z\right)
$ is a small function of $f\left( z\right) $ if $T\left( r,\varphi \right)
=S\left( r,f\right) ,$ where $S\left( r,f\right) =o\left( T\left( r,f\right)
\right) ,$ as $r\rightarrow \infty $ outside of a possible exceptional set
of finite logarithmic measure, we use $S\left( f\right) $ to denote the
family of all small functions with respect to $f\left( z\right) $. For a
meromorphic function $f\left( z\right) ,$ we define its shift by $%
f_{c}\left( z\right) =f\left( z+c\right) $ and its difference operators by%
\begin{equation*}
\Delta _{c}f\left( z\right) =f\left( z+c\right) -f\left( z\right) ,\text{ \
\ }\Delta _{c}^{n}f\left( z\right) =\Delta _{c}^{n-1}\left( \Delta
_{c}f\left( z\right) \right) ,\text{ }n\in
\mathbb{N}
,\text{ }n\geq 2.
\end{equation*}%
In particular, $\Delta _{c}^{n}f\left( z\right) =\Delta ^{n}f\left( z\right)
$ for the case $c=1.$

\quad

\noindent \qquad Let $f$ and $g$ be two meromorphic functions and let $a$ be
a finite nonzero value. We say that $f$ and $g$ share the value $a$ CM
provided that $f-a$ and $g-a$ have the same zeros counting multiplicities.
Similarly, we say that $f$ and $g$ share $a$ IM provided that $f-a$ and $g-a$
have the same zeros ignoring multiplicities. It is well-known that if $f$
and $g$ share four distinct values CM, then $f$ is a M\"{o}bius
transformation of $g$. In $\left[ 15\right] ,$ Rubel and Yang proved that if
an entire function $f$ shares two distinct complex numbers CM with its
derivative $f^{\prime }$, then $f\equiv f^{\prime }$ . In 1986, Jank et al. $%
\left( \text{see }\left[ 10\right] \right) $ proved that for a nonconstant
meromorphic function $f$, if $f$, $f^{\prime }$ and $f^{\prime \prime }$
share a finite nonzero value CM, then $f^{\prime }\equiv f$ . This result
suggests the following question:

\quad

\noindent \textbf{Question 1 }$\left[ 17\right] $\textit{\ Let }$f$\textit{\
be a nonconstant meromorphic function, let }$a$\textit{\ be a finite nonzero
constant, and let }$n$\textit{\ and }$m$\textit{\ }$\left( n<m\right) $%
\textit{\ be positive integers. If }$f$\textit{, }$f^{\left( n\right) }$%
\textit{\ and }$f^{\left( m\right) }$\textit{\ share }$a$\textit{\ CM, then
can we get the result }$f^{\left( n\right) }\equiv f$\textit{?}

\quad

\noindent The following example $\left( \text{see }\left[ 18\right] \right) $
shows that the answer to the above question is, in general, negative. Let $n$
and $m$ be positive integers satisfying $m>n+1$, and let $b$ be a constant
satisfying $b^{n}=b^{m}\neq 1$. Set $a=b^{n}$ and $f\left( z\right)
=e^{bz}+a-1$. Then $f$, $f^{\left( n\right) }$ and $f^{\left( m\right) }$
share the value $a$ CM, and $f^{\left( n\right) }\not\equiv f$. However,
when $f$ is an entire function of finite order and $m=n+1$, the answer to
Question 1 is still positive. In fact, P. Li and C. C. Yang proved the
following:

\quad

\noindent \textbf{Theorem A} $\left[ 14\right] $ \textit{Let }$f$\textit{\
be a nonconstant entire function, let }$a$\textit{\ be a finite nonzero
constant, and let }$n$\textit{\ be a positive integer. If }$f$\textit{, }$%
f^{\left( n\right) }$\textit{\ and }$f^{\left( n+1\right) }$\textit{\ share
the value }$a$\textit{\ CM, then }$f\equiv f^{\prime }.$

\textit{\quad }

\noindent \qquad Recently several papers have focussed on the Nevanlinna
theory with respect to difference operators see, e.g. $[1],$ $[5],\left[ 7%
\right] ,\left[ 8\right] $. Many authors started to investigate the
uniqueness of meromorphic functions sharing values with their shifts or
difference operators. In $\left[ 2,3\right] ,$ B. Chen et al. proved a
difference analogue of result of Jank et al. and obtained the following
results:

\quad

\noindent \textbf{Theorem B }$\left[ 2\right] $ \textit{Let }$f\left(
z\right) $ \textit{be a nonconstant entire function of finite order, and let
}$a\left( z\right) \left( \not\equiv 0\right) \in S\left( f\right) $\textit{%
\ be a periodic entire function with period }$c$\textit{. If }$f\left(
z\right) ,$\textit{\ }$\Delta _{c}f\left( z\right) $\textit{\ and }$\Delta
_{c}^{2}f\left( z\right) $\textit{\ share }$a\left( z\right) $\textit{\ CM,
then }$\Delta _{c}f\equiv \Delta _{c}^{2}f.$

\quad

\noindent \textbf{Theorem C }$\left[ 3\right] $ \textit{Let }$f\left(
z\right) $ \textit{be a nonconstant entire function of finite order, and let
}$a\left( z\right) \left( \not\equiv 0\right) \in S\left( f\right) $\textit{%
\ be a periodic entire function with period }$c$\textit{. If }$f\left(
z\right) ,$\textit{\ }$\Delta _{c}f\left( z\right) $\textit{\ and }$\Delta
_{c}^{n}f\left( z\right) $\textit{\ }$\left( n\geq 2\right) $ \textit{share }%
$a\left( z\right) $\textit{\ CM, then }$\Delta _{c}f\equiv \Delta _{c}^{n}f.$

\quad

\noindent \textbf{Theorem D }$\left[ 3\right] $ \textit{Let }$f\left(
z\right) $ \textit{be a nonconstant entire function of finite order. If }$%
f\left( z\right) ,$\textit{\ }$\Delta _{c}f\left( z\right) $\textit{\ and }$%
\Delta _{c}^{n}f\left( z\right) $\textit{\ share }$0$\textit{\ CM, then }$%
\Delta _{c}^{n}f\left( z\right) =C\Delta _{c}f\left( z\right) ,$\textit{\
where }$C$\textit{\ is a nonzero constant.}

\quad

\noindent \qquad Recently in $\left[ 13\right] ,$ Z. Latreuch et al. proved
the following results:

\quad

\noindent \textbf{Theorem E }$\left[ 13\right] $\textbf{\ }\textit{Let }$%
f\left( z\right) $ \textit{be a nonconstant entire function of finite order,
and let }$a\left( z\right) \left( \not\equiv 0\right) \in S\left( f\right) $%
\textit{\ be a periodic entire function with period }$c$\textit{. If }$%
f\left( z\right) $\textit{, }$\Delta _{c}^{n}f\left( z\right) $\textit{\ and
}$\Delta _{c}^{n+1}f\left( z\right) $\textit{\ }$\left( n\geq 1\right) $
\textit{share }$a\left( z\right) $ \textit{CM, then }$\Delta
_{c}^{n+1}f\left( z\right) \equiv \Delta _{c}^{n}f\left( z\right) .$

\quad

\noindent \textbf{Theorem F }$\left[ 13\right] $ \textit{Let }$f\left(
z\right) $ \textit{be a nonconstant entire function of finite order. If }$%
f\left( z\right) ,$\textit{\ }$\Delta _{c}^{n}f\left( z\right) $\textit{\
and }$\Delta _{c}^{n+1}f\left( z\right) $\textit{\ share }$0$\textit{\ CM,
then }$\Delta _{c}^{n+1}f\left( z\right) =C\Delta _{c}^{n}f\left( z\right) ,$%
\textit{\ where }$C$\textit{\ is a nonzero constant.}

\quad

\noindent \qquad For the case $n=1,$ A. El Farissi and others gave the
following improvement.

\quad

\noindent \textbf{Theorem G }$\left[ 6\right] $ \textit{Let }$f\left(
z\right) $ \textit{be a non-periodic entire function of finite order, and
let }$a\left( z\right) \left( \not\equiv 0\right) \in S\left( f\right) $%
\textit{\ be a periodic entire function with period }$c$\textit{. If }$%
f\left( z\right) ,$\textit{\ }$\Delta _{c}f\left( z\right) $\textit{\ and }$%
\Delta _{c}^{2}f\left( z\right) $\textit{\ share }$a\left( z\right) $\textit{%
\ CM, then }$\Delta _{c}f\left( z\right) \equiv f\left( z\right) .$

\quad

\noindent \textbf{Remark 1.1} The Theorem G is essentially known in [6]. For
the convenience of readers, we give his proof in the Lemma 2.4.

\quad

\noindent \qquad It is naturally now to ask the following question: \textit{%
Under the hypotheses of Theorem E, can we get the result }$\Delta
_{c}f\left( z\right) \equiv f\left( z\right) $\textit{? }The aim of this
paper is to answer this question and to give a difference analogue of result
of P. Li and C. C. Yang in $\left[ 14\right] $. In fact we obtain the
following results:

\quad

\noindent \textbf{Theorem 1.1} \textbf{\ }\textit{Let }$f\left( z\right) $
\textit{be a nonconstant entire function of finite order such that }$\Delta
_{c}^{n}f\left( z\right) \not\equiv 0$\textit{, and let }$a\left( z\right)
\left( \not\equiv 0\right) \in S\left( f\right) $\textit{\ be a periodic
entire function with period }$c$\textit{. If }$f\left( z\right) $\textit{, }$%
\Delta _{c}^{n}f\left( z\right) $\textit{\ and }$\Delta _{c}^{n+1}f\left(
z\right) $\textit{\ }$\left( n\geq 1\right) $ \textit{share }$a\left(
z\right) $ \textit{CM, then }$\Delta _{c}f\left( z\right) \equiv f\left(
z\right) .$

\quad

\noindent \textbf{Remark 1.2 }The condition\textbf{\ }$\Delta
_{c}^{n}f\left( z\right) \not\equiv 0$ is necessary. Let's take for example
the entire function $f\left( z\right) =1+e^{2\pi iz}$ and $c=a=1,$ then $f-a$
and $\Delta ^{n}f-a=\Delta ^{n+1}f-a=-1$ have the same zeros but $\Delta
f\neq f.$ On the other hand, under the conditions of Theorem 1.1 $\Delta
_{c}^{n}f\left( z\right) \not\equiv 0$ can not be a periodic entire function
because $\Delta _{c}^{n+1}f\left( z\right) \equiv \Delta _{c}^{n}f\left(
z\right) $ (Theorem E, $\left[ 13\right] $).

\quad

\noindent \textbf{Example 1.1 }Let $f\left( z\right) =e^{z\ln 2}$ and $c=1.$
Then, for any $a\in
\mathbb{C}
,$ we notice that $f\left( z\right) ,$ $\Delta _{c}^{n}f\left( z\right) $%
\textit{\ }and $\Delta _{c}^{n+1}f\left( z\right) $\ share $a$\ CM for all $%
n\in
\mathbb{N}
$ and we can easily see that $\Delta _{c}f\left( z\right) \equiv f\left(
z\right) .$ This example satisfies Theorem 1.1.

\quad

\noindent \textbf{Theorem 1.2 }\textit{Let }$f\left( z\right) $ \textit{be a
nonconstant entire function of finite order such that }$\Delta
_{c}^{n}f\left( z\right) \not\equiv 0$\textit{, and let }$a\left( z\right) ,$
$b\left( z\right) \left( \not\equiv 0\right) \in S\left( f\right) $\textit{\
such that }$b\left( z\right) $ \textit{is a periodic entire function with
period }$c$ \textit{and }$\Delta _{c}^{m}a\left( z\right) \equiv 0$ $\left(
1\leq m\leq n\right) $\textit{. If }$f\left( z\right) -a\left( z\right) ,$%
\textit{\ }$\Delta _{c}^{n}f\left( z\right) -b\left( z\right) $\textit{\ and
}$\Delta _{c}^{n+1}f\left( z\right) -b\left( z\right) $\textit{\ share }$0$%
\textit{\ CM, then }$\Delta _{c}f\left( z\right) \equiv f\left( z\right)
+b\left( z\right) +\Delta _{c}a\left( z\right) -a\left( z\right) .$

\quad

\noindent \textbf{Remark 1.3 }The condition $b\left( z\right) \not\equiv 0$
is necessary in the proof of Theorem 1.2, for the case $b\left( z\right)
\equiv 0,$ please see Theorem 1.4.

\quad

\noindent \textbf{Remark 1.4 }The condition $\Delta _{c}^{m}a\left( z\right)
\equiv 0$ in Theorem 1.2 is more general than the condition "periodic entire
function of period $c$".

\quad

\noindent For the case $m=1,$ we deduce the following result.

\quad

\noindent \textbf{Corollary 1.1 }\textit{Let }$f\left( z\right) $ \textit{be
a nonconstant entire function of finite order such that }$\Delta
_{c}^{n}f\left( z\right) \not\equiv 0$\textit{, and let }$a\left( z\right) ,$
$b\left( z\right) \left( \not\equiv 0\right) \in S\left( f\right) $\textit{\
be periodic entire functions with period }$c$\textit{. If }$f\left( z\right)
-a\left( z\right) ,$\textit{\ }$\Delta _{c}^{n}f\left( z\right) -b\left(
z\right) $\textit{\ and }$\Delta _{c}^{n+1}f\left( z\right) -b\left(
z\right) $\textit{\ share }$0$\textit{\ CM, then }$\Delta _{c}f\left(
z\right) \equiv f\left( z\right) +b\left( z\right) -a\left( z\right) .$

\quad

\noindent \textbf{Example 1.2 }Let $f\left( z\right) =e^{z\ln 2}-2,$ $a=-1$
and $b=1.$ It is clear that $f\left( z\right) -a,$ $\Delta ^{n}f\left(
z\right) -b$ and $\Delta ^{n+1}f\left( z\right) -b$ share $0$ CM. Here, we
also get $\Delta f\left( z\right) =f\left( z\right) +b-a.$

\quad

\noindent \textbf{Example 1.3 }Let $f\left( z\right) =e^{z\ln 2}+z^{3}-1,$ $%
a\left( z\right) =z^{3}$ and $b=1.$ It is clear that $f\left( z\right)
-z^{3},$ $\Delta ^{4}f\left( z\right) -1$ and $\Delta ^{5}f\left( z\right)
-1 $ share $0$ CM. On the other hand, we can verify that $\Delta f\left(
z\right) =f\left( z\right) +1+\Delta z^{3}-z^{3}$ which satisfies Theorem
1.2.

\quad

\noindent \textbf{Theorem 1.3 }\textit{Let }$f\left( z\right) $ \textit{be a
nonconstant entire function of finite order such that }$\Delta
_{c}^{n}f\left( z\right) \not\equiv 0$\textit{. If }$f\left( z\right) ,$%
\textit{\ }$\Delta _{c}^{n}f\left( z\right) $\textit{\ and }$\Delta
_{c}^{n+1}f\left( z\right) $\textit{\ share }$0$\textit{\ CM, then }$\Delta
_{c}f\left( z\right) \equiv Cf\left( z\right) ,$\textit{\ where }$C$\textit{%
\ is a nonzero constant.}

\quad

\noindent \textbf{Example 1.4 }Let $f\left( z\right) =e^{az}$ and $c=1$
where $a\neq 2k\pi i$ $\left( k\in
\mathbb{Z}
\right) ,$ it is clear that $\Delta _{c}^{n}f\left( z\right) =\left(
e^{a}-1\right) ^{n}e^{az}$ for any integer $n\geq 1.$ So, $f\left( z\right)
, $ $\Delta _{c}^{n}f\left( z\right) $\textit{\ }and $\Delta
_{c}^{n+1}f\left( z\right) $\ share $0$\ CM for all $n\in
\mathbb{N}
$ and we can easily see that $\Delta _{c}f\left( z\right) \equiv Cf\left(
z\right) $ where $C=e^{a}-1.$ This example satisfies Theorem 1.3.

\quad

\noindent \textbf{Corollary 1.2 }\textit{Let }$f\left( z\right) $ \textit{be
a nonconstant entire function of finite order such that }$f\left( z\right) ,$%
\textit{\ }$\Delta _{c}^{n}f\left( z\right) $\textit{\ }$\left( \not\equiv
0\right) $ \textit{and }$\Delta _{c}^{n+1}f\left( z\right) $\textit{\ }$%
\left( n\geq 1\right) $ \textit{share }$0$\textit{\ CM. If there exists a
point }$z_{0}$ \textit{and} \textit{an integer} $m\geq 1$ \textit{such that }%
$\Delta _{c}^{m}f\left( z_{0}\right) =f\left( z_{0}\right) \neq 0,$ \textit{%
then }$\Delta _{c}^{m}f\left( z\right) \equiv f\left( z\right) .$

\quad

\noindent By combining Theorem 1.2 and Theorem 1.3 we can prove the
following result.

\quad

\noindent \textbf{Theorem 1.4 }\textit{Let }$f\left( z\right) $ \textit{be a
nonconstant entire function of finite order such that }$\Delta
_{c}^{n}f\left( z\right) \not\equiv 0$\textit{, and let }$a\left( z\right)
\in S\left( f\right) $\textit{\ such that }$\Delta _{c}^{m}a\left( z\right)
\equiv 0$ $\left( 1\leq m\leq n\right) $\textit{. If }$f\left( z\right)
-a\left( z\right) ,$\textit{\ }$\Delta _{c}^{n}f\left( z\right) $\textit{\
and }$\Delta _{c}^{n+1}f\left( z\right) $\textit{\ share }$0$\textit{\ CM,
then }$\Delta _{c}f\left( z\right) \equiv Cf\left( z\right) +\Delta
_{c}a\left( z\right) -a\left( z\right) ,$\textit{\ where }$C$\textit{\ is a
nonzero constant.}

\section{Some lemmas}

\noindent \textbf{Lemma 2.1 }$\left[ 5\right] $ \textit{Let }$\eta _{1},\eta
_{2}$ \textit{be two arbitrary complex numbers such that }$\eta _{1}\neq
\eta _{2}$ \textit{and let }$f\left( z\right) $ \textit{be a finite order
meromorphic function. Let }$\sigma $ \textit{be the order of }$f\left(
z\right) $, \textit{then for each }$\varepsilon >0,$ \textit{we have}%
\begin{equation*}
m\left( r,\frac{f\left( z+\eta _{1}\right) }{f\left( z+\eta _{2}\right) }%
\right) =O\left( r^{\sigma -1+\varepsilon }\right) .
\end{equation*}

\noindent \qquad By combining Theorem 1.4 in $\left[ 4\right] $ and Theorem
2.2 in $\left[ 12\right] ,$ we can prove the following lemma.

\quad

\noindent \textbf{Lemma 2.2 }\textit{Let }$a_{0}\left( z\right) ,a_{1}\left(
z\right) ,\cdots ,a_{n}\left( z\right) \left( \not\equiv 0\right) ,$\textit{%
\ }$F\left( z\right) \left( \not\equiv 0\right) $\textit{\ be finite order
meromorphic functions, }$c_{k}$\textit{\ (}$k=0,\cdots ,n$\textit{) be
constants, unequal to each other. If }$f$\textit{\ is a finite order
meromorphic solution of the equation}%
\begin{equation}
a_{n}\left( z\right) f\left( z+c_{n}\right) +\cdots +a_{1}\left( z\right)
f\left( z+c_{1}\right) +a_{0}\left( z\right) f\left( z+c_{0}\right) =F\left(
z\right)  \tag{2.1}
\end{equation}%
\textit{with }%
\begin{equation*}
\max \left\{ \rho \left( a_{i}\right) ,\left( i=0,\cdots ,n\right) ,\rho
\left( F\right) \right\} <\rho \left( f\right) ,
\end{equation*}%
\textit{then }$\lambda \left( f\right) =\rho \left( f\right) .$

\quad

\noindent \textit{Proof.}\textbf{\ }By $\left( 2.1\right) $ we have
\begin{equation}
\frac{1}{f\left( z+c_{0}\right) }=\frac{1}{F}\left( a_{n}\frac{f\left(
z+c_{n}\right) }{f\left( z+c_{0}\right) }+\cdots +a_{1}\frac{f\left(
z+c_{1}\right) }{f\left( z+c_{0}\right) }+a_{0}\right) .  \tag{2.2}
\end{equation}%
Set $\max \left\{ \rho \left( a_{j}\right) \text{ }\left( j=0,\cdots
,n\right) ,\rho \left( F\right) \right\} =\beta <\rho \left( f\right) =\rho
. $ Then, for any given $\varepsilon $ $\left( 0<\varepsilon <\frac{\rho
-\beta }{2}\right) $, we have%
\begin{equation}
\overset{n}{\underset{j=0}{\sum }}T\left( r,a_{j}\right) +T\left( r,F\right)
\leq \left( n+2\right) \exp \left\{ r^{\beta +\varepsilon }\right\} =o\left(
1\right) \exp \left\{ r^{\rho -\varepsilon }\right\} .  \tag{2.3}
\end{equation}%
By $\left( 2.2\right) $, $\left( 2.3\right) $ and Lemma 2.1, we obtain%
\begin{equation*}
T\left( r,f\right) =T\left( r,\frac{1}{f}\right) +O\left( 1\right) =m\left(
r,\frac{1}{f}\right) +N\left( r,\frac{1}{f}\right) +O\left( 1\right)
\end{equation*}%
\begin{equation*}
\leq N\left( r,\frac{1}{f}\right) +m\left( r,\frac{1}{F}\right) +\overset{n}{%
\underset{j=0}{\sum }}m\left( r,a_{j}\right) +\overset{n}{\underset{j=1}{%
\sum }}m\left( r,\frac{f\left( z+c_{j}\right) }{f\left( z+c_{0}\right) }%
\right) +O\left( 1\right)
\end{equation*}%
\begin{equation*}
\leq N\left( r,\frac{1}{f}\right) +T\left( r,\frac{1}{F}\right) +\overset{n}{%
\underset{j=0}{\sum }}T\left( r,a_{j}\right) +\overset{n}{\underset{j=1}{%
\sum }}m\left( r,\frac{f\left( z+c_{j}\right) }{f\left( z+c_{0}\right) }%
\right) +O\left( 1\right)
\end{equation*}%
\begin{equation}
\leq N\left( r,\frac{1}{f}\right) +O\left( r^{\rho -1+\varepsilon }\right)
+o\left( 1\right) \exp \left\{ r^{\rho -\varepsilon }\right\} .  \tag{2.4}
\end{equation}%
By $\left( 2.4\right) $, we obtain that $\rho \left( f\right) \leq \lambda
\left( f\right) $ and since $\lambda \left( f\right) \leq \rho \left(
f\right) $ for every meromorphic function, we deduce that $\lambda \left(
f\right) =\rho \left( f\right) .$

\quad

\noindent \textbf{Remark 2.1 }Recently, Shun-Zhou Wu and Xiu-Min Zheng (see $%
\left[ 16\right] $) obtained Lemma 2.2 by using a different proof.

\quad

\noindent \textbf{Lemma 2.3 }$\left[ 17\right] $ \textit{Suppose }$%
f_{j}\left( z\right) $\textit{\ }$(j=1,2,\cdots ,n+1)$\textit{\ and }$%
g_{j}\left( z\right) $\textit{\ }$(j=1,2,\cdots ,n)$\textit{\ }$(n\geq 1)$%
\textit{\ are entire functions satisfying the following conditions:}

\noindent $\left( \text{i}\right) $\textit{\ }$\overset{n}{\underset{j=1}{%
\sum }}f_{j}\left( z\right) e^{g_{j}\left( z\right) }\equiv f_{n+1}\left(
z\right) ;$

\noindent $\left( \text{ii}\right) $\textit{\ The order of }$f_{j}\left(
z\right) $\textit{\ is less than the order of }$e^{g_{k}\left( z\right) }$%
\textit{\ for }$1\leq j\leq n+1,$\textit{\ }$1\leq k\leq n.$\textit{\ And
furthermore, the order of }$f_{j}\left( z\right) $\textit{\ is less than the
order of }$e^{g_{h}\left( z\right) -g_{k}\left( z\right) }$\textit{\ for }$%
n\geq 2$\textit{\ and }$1\leq j\leq n+1,$\textit{\ }$1\leq h<k\leq n.$

\noindent \textit{Then }$f_{j}\left( z\right) \equiv 0,$\textit{\ }$\left(
j=1,2,\cdots n+1\right) .$

\quad

\noindent \textbf{Lemma 2.4 }$\left[ 6\right] $ \textit{Let }$f\left(
z\right) $ \textit{be a non-periodic entire function of finite order, and
let }$a\left( z\right) \left( \not\equiv 0\right) \in S\left( f\right) $%
\textit{\ be a periodic entire function with period }$c$\textit{. If }$%
f\left( z\right) ,$\textit{\ }$\Delta _{c}f\left( z\right) $\textit{\ and }$%
\Delta _{c}^{2}f\left( z\right) $\textit{\ share }$a\left( z\right) $\textit{%
\ CM, then }$\Delta _{c}f\left( z\right) \equiv f\left( z\right) .$

\quad

\noindent \textit{Proof.} Suppose that $\Delta _{c}f\left( z\right)
\not\equiv f\left( z\right) .$ Since $f,$ $\Delta _{c}f$ and $\Delta
_{c}^{2}f$ share $a\left( z\right) $ CM, then we have%
\begin{equation*}
\frac{\Delta _{c}f\left( z\right) -a\left( z\right) }{f\left( z\right)
-a\left( z\right) }=e^{P\left( z\right) }
\end{equation*}%
and
\begin{equation*}
\frac{\Delta _{c}^{2}f\left( z\right) -a\left( z\right) }{f\left( z\right)
-a\left( z\right) }=e^{Q\left( z\right) }
\end{equation*}%
where $P$ $\left( e^{P}\not\equiv 1\right) $ and $Q$ are polynomials. By
using Theorem B, we obtain that $\Delta _{c}^{2}f\equiv \Delta _{c}f$, which
means that
\begin{equation}
\alpha \left( z\right) =\Delta _{c}f\left( z\right) -f\left( z\right)
\tag{2.5}
\end{equation}%
is entire periodic function of period $c.$ By $\left( 2.5\right) $ we have%
\begin{equation*}
\Delta _{c}f\left( z\right) -a\left( z\right) =f\left( z\right) -a\left(
z\right) +\alpha \left( z\right) ,
\end{equation*}%
then%
\begin{equation*}
\frac{\Delta _{c}f\left( z\right) -a\left( z\right) }{f\left( z\right)
-a\left( z\right) }=1+\frac{\alpha \left( z\right) }{f\left( z\right)
-a\left( z\right) }=e^{P\left( z\right) },
\end{equation*}%
which is equivalent to%
\begin{equation}
f\left( z\right) -a\left( z\right) =\frac{\alpha \left( z\right) }{%
e^{P\left( z\right) }-1}.  \tag{2.6}
\end{equation}%
Since $\alpha \left( z\right) $ and $a\left( z\right) $ are periodic
functions of period $c$, then we have
\begin{equation}
\Delta _{c}f\left( z\right) =\alpha \left( z\right) \Delta _{c}\left( \frac{1%
}{e^{P\left( z\right) }-1}\right)  \tag{2.7}
\end{equation}%
and
\begin{equation}
\Delta _{c}^{2}f\left( z\right) =\alpha \left( z\right) \Delta
_{c}^{2}\left( \frac{1}{e^{P\left( z\right) }-1}\right) .  \tag{2.8}
\end{equation}%
We have the two following subcases:

\noindent $\left( \text{i}\right) $ If $P\equiv K$ $\left( K\neq 2k\pi i,%
\text{ }K\in
\mathbb{Z}
\right) $, then by $\left( 2.7\right) $ we have $\Delta _{c}f\left( z\right)
=0.$ On the other hand, by using $\left( 2.5\right) ,$ $\left( 2.6\right) $
and $\Delta _{c}f\left( z\right) =0,$ we deduce that
\begin{equation*}
f\left( z\right) -a\left( z\right) =\frac{-f\left( z\right) }{e^{K}-1},\text{
}K\in
\mathbb{C}
-\left\{ 2k\pi i,k\in
\mathbb{Z}
\right\} .
\end{equation*}%
So,
\begin{equation*}
f\left( z\right) =\frac{e^{K}-1}{e^{K}}a\left( z\right) .
\end{equation*}%
Hence
\begin{equation*}
T\left( r,f\right) =S\left( r,f\right) ,
\end{equation*}%
which is a contradiction.

\noindent $\left( \text{ii}\right) $ If $P$ is nonconstant and since $\Delta
_{c}^{2}f\left( z\right) =\Delta _{c}f\left( z\right) ,$ then
\begin{equation*}
e^{P_{c}\left( z\right) +P\left( z\right) }-3e^{P_{2c}\left( z\right)
+P\left( z\right) }+2e^{P_{2c}\left( z\right) +P_{c}\left( z\right)
}+e^{P_{2c}\left( z\right) }-3e^{P_{c}\left( z\right) }+2e^{P\left( z\right)
}=0
\end{equation*}%
which is equivalent to%
\begin{equation}
e^{P_{c}\left( z\right) }+\left( 2e^{\Delta _{c}P\left( z\right) }-3\right)
e^{P_{2c}\left( z\right) }=-e^{\Delta _{c}P_{c}\left( z\right) +\Delta
_{c}P\left( z\right) }+3e^{\Delta _{c}P\left( z\right) }-2.  \tag{2.9}
\end{equation}%
Since $\deg \Delta _{c}P=\deg P-1,$ then we have
\begin{equation}
\rho \left( e^{P_{c}}+\left( 2e^{\Delta _{c}P}-3\right) e^{P_{2c}}\right)
=\rho \left( -e^{\Delta _{c}P_{c}+\Delta _{c}P}+3e^{\Delta _{c}P}-2\right)
\leq \deg P-1.  \tag{2.10}
\end{equation}%
On the other hand,
\begin{equation}
\rho \left( e^{P_{c}}+\left( 2e^{\Delta _{c}P}-3\right) e^{P_{2c}}\right)
=\rho \left( e^{P_{c}}\right) =\deg P  \tag{2.11}
\end{equation}%
because if we have the contrary
\begin{equation*}
\rho \left( e^{P_{c}}+\left( 2e^{\Delta _{c}P}-3\right) e^{P_{2c}}\right)
<\rho \left( e^{P_{c}}\right) ,
\end{equation*}%
we obtain the following contradiction%
\begin{equation*}
\deg P=\rho \left( \frac{e^{P_{c}}+\left( 2e^{\Delta _{c}P}-3\right)
e^{P_{2c}}}{e^{P_{c}}}\right) =\rho \left( 1+\left( 2e^{\Delta
_{c}P}-3\right) e^{\Delta P_{c}}\right) \leq \deg P-1.
\end{equation*}%
By using $\left( 2.10\right) $ and $\left( 2.11\right) ,$ we obtain%
\begin{equation*}
\deg P\leq \deg P-1
\end{equation*}%
which is a contradiction. This leads to $\Delta _{c}f\left( z\right)
=f\left( z\right) $. Thus, the proof of Lemma 2.4 is completed.

\section{Proof of the Theorems and Corollary}

\noindent \textbf{Proof of the Theorem 1.1.} Obviously, suppose that $\Delta
_{c}f\left( z\right) \not\equiv f\left( z\right) $. By using Theorem E, we
have
\begin{equation}
\frac{\Delta _{c}^{n}f\left( z\right) -a\left( z\right) }{f\left( z\right)
-a\left( z\right) }=e^{P\left( z\right) }  \tag{3.1}
\end{equation}%
and
\begin{equation}
\frac{\Delta _{c}^{n+1}f\left( z\right) -a\left( z\right) }{f\left( z\right)
-a\left( z\right) }=e^{P\left( z\right) },  \tag{3.2}
\end{equation}%
where $P$ $\left( e^{P}\not\equiv 1\right) $ is polynomial. Dividing the
proof of Theorem 1.1 into two cases:

\quad

\noindent \textbf{Case} \textbf{1}. $P$ is a nonconstant polynomial. Setting
now $g\left( z\right) =f\left( z\right) -a\left( z\right) .$ Then, we have
from $\left( 3.1\right) $ and $\left( 3.2\right) $%
\begin{equation}
\Delta _{c}^{n}g\left( z\right) =e^{P\left( z\right) }g\left( z\right)
+a\left( z\right)  \tag{3.3}
\end{equation}%
and
\begin{equation}
\Delta _{c}^{n+1}g\left( z\right) =e^{P\left( z\right) }g\left( z\right)
+a\left( z\right) .  \tag{3.4}
\end{equation}%
By $\left( 3.3\right) $ and $\left( 3.4\right) ,$ we have%
\begin{equation*}
g_{c}\left( z\right) =2e^{P-P_{c}}g\left( z\right) +a\left( z\right)
e^{-P_{c}}.
\end{equation*}%
Using the principle of mathematical induction, we obtain
\begin{equation}
g_{ic}\left( z\right) =2^{i}e^{P-P_{ic}}g\left( z\right) +a\left( z\right)
\left( 2^{i}-1\right) e^{-P_{ic}},\text{ }i\geq 1.  \tag{3.5}
\end{equation}%
Now, we can rewrite $\left( 3.3\right) $ as
\begin{equation*}
\Delta _{c}^{n}g\left( z\right) =\overset{n}{\underset{i=1}{\sum }}%
C_{n}^{i}\left( -1\right) ^{n-i}\left( 2^{i}e^{P-P_{ic}}g\left( z\right)
+a\left( z\right) \left( 2^{i}-1\right) e^{-P_{ic}}\right)
\end{equation*}%
\begin{equation*}
+\left( -1\right) ^{n}g\left( z\right) =e^{P}g\left( z\right) +a\left(
z\right) ,
\end{equation*}%
which implies%
\begin{equation*}
\left( \overset{n}{\underset{i=0}{\sum }}C_{n}^{i}\left( -1\right)
^{n-i}2^{i}e^{P-P_{ic}}-e^{P}\right) g\left( z\right)
\end{equation*}%
\begin{equation*}
+a\left( z\right) \left( \overset{n}{\underset{i=0}{\sum }}C_{n}^{i}\left(
-1\right) ^{n-i}\left( 2^{i}-1\right) e^{-P_{ic}}-1\right) =0.
\end{equation*}%
Hence
\begin{equation}
A_{n}\left( z\right) g\left( z\right) +B_{n}\left( z\right) =0,  \tag{3.6}
\end{equation}%
where
\begin{equation*}
A_{n}\left( z\right) =\overset{n}{\underset{i=0}{\sum }}C_{n}^{i}\left(
-1\right) ^{n-i}2^{i}e^{P-P_{ic}}-e^{P}
\end{equation*}%
and
\begin{equation*}
B_{n}\left( z\right) =a\left( z\right) \left( \overset{n}{\underset{i=0}{%
\sum }}C_{n}^{i}\left( -1\right) ^{n-i}\left( 2^{i}-1\right)
e^{-P_{ic}}-1\right) .
\end{equation*}%
By the same method, we can rewrite $\left( 3.4\right) $ as%
\begin{equation}
A_{n+1}\left( z\right) g\left( z\right) +B_{n+1}\left( z\right) =0,
\tag{3.7}
\end{equation}%
where
\begin{equation*}
A_{n+1}\left( z\right) =\overset{n+1}{\underset{i=0}{\sum }}%
C_{n+1}^{i}\left( -1\right) ^{n+1-i}2^{i}e^{P-P_{ic}}-e^{P}
\end{equation*}%
and
\begin{equation*}
B_{n+1}\left( z\right) =a\left( z\right) \left( \overset{n+1}{\underset{i=0}{%
\sum }}C_{n+1}^{i}\left( -1\right) ^{n+1-i}\left( 2^{i}-1\right)
e^{-P_{ic}}-1\right) .
\end{equation*}%
We can see easily from the equations $\left( 3.6\right) $ and $\left(
3.7\right) $ that
\begin{equation}
h\left( z\right) =A_{n}\left( z\right) B_{n+1}\left( z\right) -A_{n+1}\left(
z\right) B_{n}\left( z\right) \equiv 0.  \tag{3.8}
\end{equation}%
On the other hand, we remark that
\begin{equation*}
e^{P}B_{n}\left( z\right) =a\left( z\right) e^{P}\left( \overset{n}{\underset%
{i=0}{\sum }}C_{n}^{i}\left( -1\right) ^{n-i}2^{i}e^{-P_{ic}}-\overset{n}{%
\underset{i=0}{\sum }}C_{n}^{i}\left( -1\right) ^{n-i}e^{-P_{ic}}-1\right)
\end{equation*}%
\begin{equation*}
=a\left( z\right) e^{P}\left( \overset{n}{\underset{i=0}{\sum }}%
C_{n}^{i}\left( -1\right) ^{n-i}2^{i}e^{-P_{ic}}-1-\Delta _{c}^{n}\left(
e^{-P}\right) \right)
\end{equation*}%
\begin{equation*}
=a\left( z\right) \left( A_{n}\left( z\right) -e^{P}\Delta _{c}^{n}\left(
e^{-P}\right) \right) .
\end{equation*}%
Then
\begin{equation}
B_{n}\left( z\right) =a\left( z\right) \left( e^{-P}A_{n}\left( z\right)
-\Delta _{c}^{n}\left( e^{-P}\right) \right) .  \tag{3.9}
\end{equation}%
By the same method, we obtain%
\begin{equation}
B_{n+1}\left( z\right) =a\left( z\right) \left( e^{-P}A_{n+1}\left( z\right)
-\Delta _{c}^{n+1}\left( e^{-P}\right) \right) .  \tag{3.10}
\end{equation}%
Return now to the equation $\left( 3.8\right) $, by using $\left( 3.9\right)
$ and $\left( 3.10\right) ,$ we get%
\begin{equation*}
h\left( z\right) =A_{n}\left( z\right) B_{n+1}\left( z\right) -A_{n+1}\left(
z\right) B_{n}\left( z\right)
\end{equation*}%
\begin{equation*}
=A_{n}\left( z\right) \left[ a\left( z\right) \left( e^{-P}A_{n+1}\left(
z\right) -\Delta _{c}^{n+1}\left( e^{-P}\right) \right) \right]
\end{equation*}%
\begin{equation*}
-A_{n+1}\left( z\right) \left[ a\left( z\right) \left( e^{-P}A_{n}\left(
z\right) -\Delta _{c}^{n}\left( e^{-P}\right) \right) \right]
\end{equation*}%
\begin{equation*}
=a\left( z\right) \left[ A_{n+1}\left( z\right) \Delta _{c}^{n}\left(
e^{-P}\right) -A_{n}\left( z\right) \Delta _{c}^{n+1}\left( e^{-P}\right) %
\right] \equiv 0.
\end{equation*}%
Hence%
\begin{equation*}
A_{n+1}\left( z\right) \Delta _{c}^{n}\left( e^{-P}\right) -A_{n}\left(
z\right) \Delta _{c}^{n+1}\left( e^{-P}\right) \equiv 0.
\end{equation*}%
Therefore%
\begin{equation*}
\Delta _{c}^{n}\left( e^{-P}\right) \left( \overset{n+1}{\underset{i=0}{\sum
}}C_{n+1}^{i}\left( -1\right) ^{n+1-i}2^{i}e^{-P_{ic}}-1\right)
\end{equation*}%
\begin{equation*}
-\Delta _{c}^{n+1}\left( e^{-P}\right) \left( \overset{n}{\underset{i=0}{%
\sum }}C_{n}^{i}\left( -1\right) ^{n-i}2^{i}e^{-P_{ic}}-1\right) =0.
\end{equation*}%
Thus%
\begin{eqnarray*}
&&\Delta _{c}^{n}\left( e^{-P}\right) \overset{n+1}{\underset{i=0}{\sum }}%
C_{n+1}^{i}\left( -1\right) ^{n+1-i}2^{i}e^{-P_{i}c}-\Delta _{c}^{n+1}\left(
e^{-P}\right) \overset{n}{\underset{i=0}{\sum }}C_{n}^{i}\left( -1\right)
^{n-i}2^{i}e^{-P_{ic}} \\
&=&\Delta _{c}^{n}\left( e^{-P}\right) -\Delta _{c}^{n+1}\left(
e^{-P}\right) =\Delta _{c}^{n}\left( 2e^{-P}-e^{-P_{c}}\right) .
\end{eqnarray*}%
Then%
\begin{equation*}
\overset{n}{\underset{i=0}{\sum }}\left( \Delta _{c}^{n}\left( e^{-P}\right)
C_{n+1}^{i}\left( -1\right) ^{n+1-i}-\Delta _{c}^{n+1}\left( e^{-P}\right)
C_{n}^{i}\left( -1\right) ^{n-i}\right) 2^{i}e^{-P_{ic}}
\end{equation*}%
\begin{equation*}
+\Delta _{c}^{n}\left( e^{-P}\right) 2^{n+1}e^{-P_{\left( n+1\right)
c}}=\Delta _{c}^{n}\left( 2e^{-P}-e^{-P_{c}}\right) ,
\end{equation*}%
which yields%
\begin{equation*}
\overset{n}{\underset{i=0}{\sum }}\left( \Delta _{c}^{n}\left( e^{-P}\right)
C_{n+1}^{i}+\Delta _{c}^{n+1}\left( e^{-P}\right) C_{n}^{i}\right) \left(
-1\right) ^{n+1-i}2^{i}e^{P_{\left( n+1\right) c}-P_{ic}}
\end{equation*}%
\begin{equation}
+\Delta _{c}^{n}\left( e^{-P}\right) 2^{n+1}=e^{P_{\left( n+1\right)
c}}\Delta _{c}^{n}\left( 2e^{-P}-e^{-P_{c}}\right) .  \tag{3.11}
\end{equation}%
Let us denote
\begin{equation*}
\alpha _{i}\left( z\right) =\left( -1\right) ^{n+1-i}2^{i}e^{P_{\left(
n+1\right) c}-P_{ic}},\text{ }i=0,\cdots ,n
\end{equation*}%
and
\begin{equation*}
\alpha _{n+1}\left( z\right) =e^{P_{\left( n+1\right) c}}\Delta
_{c}^{n}\left( 2e^{-P}-e^{-P_{c}}\right) .
\end{equation*}%
It is clear that $\rho \left( \alpha _{i}\right) \leq \deg P-1$ for all $%
i=0,2,\cdots ,n+1.$ The equation $\left( 3.11\right) $ will be%
\begin{equation*}
\overset{n}{\underset{i=0}{\sum }}\left( \Delta _{c}^{n}\left( e^{-P}\right)
C_{n+1}^{i}+\Delta _{c}^{n+1}\left( e^{-P}\right) C_{n}^{i}\right) \alpha
_{i}\left( z\right) +\Delta _{c}^{n}\left( e^{-P}\right) 2^{n+1}
\end{equation*}%
\begin{equation*}
=\left( \overset{n}{\underset{i=0}{\sum }}C_{n+1}^{i}\alpha _{i}\left(
z\right) +2^{n+1}\right) \Delta _{c}^{n}\left( e^{-P}\right)
\end{equation*}%
\begin{equation}
+\left( \overset{n}{\underset{i=0}{\sum }}C_{n}^{i}\alpha _{i}\left(
z\right) \right) \Delta _{c}^{n+1}\left( e^{-P}\right) =\alpha _{n+1}\left(
z\right) .  \tag{3.12}
\end{equation}%
For convenience, we denote by $M\left( z\right) $ and $N\left( z\right) $
the following
\begin{equation*}
M\left( z\right) =\overset{n}{\underset{i=0}{\sum }}C_{n+1}^{i}\alpha
_{i}\left( z\right) +2^{n+1},\text{ }N\left( z\right) =\overset{n}{\underset{%
i=0}{\sum }}C_{n}^{i}\alpha _{i}\left( z\right) .
\end{equation*}%
Equation $\left( 3.12\right) $ is equivalent to
\begin{equation*}
M\left( z\right) \overset{n}{\underset{i=0}{\sum }}C_{n}^{i}\left( -1\right)
^{n-i}e^{-P_{ic}}+N\left( z\right) \overset{n+1}{\underset{i=0}{\sum }}%
C_{n+1}^{i}\left( -1\right) ^{n+1-i}e^{-P_{ic}}
\end{equation*}%
\begin{equation}
=\overset{n}{\underset{i=0}{\sum }}\left( C_{n}^{i}M\left( z\right)
-C_{n+1}^{i}N\left( z\right) \right) \left( -1\right)
^{n-i}e^{-P_{ic}}+N\left( z\right) e^{-P_{\left( n+1\right) c}}=\alpha
_{n+1}\left( z\right) .  \tag{3.13}
\end{equation}%
As conclusion, we can say that $\left( 3.13\right) $ can be written as follow%
\begin{equation}
a_{n+1}\left( z\right) e^{-P\left( z+\left( n+1\right) c\right)
}+a_{n}\left( z\right) e^{-P\left( z+nc\right) }+\cdots +a_{0}\left(
z\right) e^{-P\left( z\right) }=\alpha _{n+1}\left( z\right) ,  \tag{3.14}
\end{equation}%
where $a_{0}\left( z\right) ,\cdots ,a_{n+1}\left( z\right) $ and $\alpha
_{n+1}\left( z\right) $ are entire functions. We distingue the following two
subcases.

\quad

\noindent $\left( \text{i}\right) $ If $\deg P>1,$ then we have
\begin{equation}
\max \left\{ \rho \left( a_{i}\right) \text{ }\left( i=0,\cdots ,n+1\right)
,\rho \left( \alpha _{n+1}\right) \right\} <\deg P.  \tag{3.15}
\end{equation}%
In order to prove that $\alpha _{n+1}\left( z\right) \not\equiv 0,$ it
suffices to show that $\Delta _{c}^{n}\left( 2e^{-P}-e^{-P_{c}}\right)
\not\equiv 0.$ Suppose the contrary. Thus%
\begin{equation}
\overset{n}{\underset{i=0}{\sum }}C_{n}^{i}\left( -1\right) ^{n-i}\left(
2e^{-P_{ic}}-e^{-P_{\left( i+1\right) c}}\right) \equiv 0.  \tag{3.16}
\end{equation}%
The equation $\left( 3.16\right) $ can be written as
\begin{equation*}
\overset{n+1}{\underset{i=0}{\sum }}b_{i}e^{-P_{ic}}\equiv 0,
\end{equation*}%
where
\begin{equation*}
b_{i}=\left\{
\begin{array}{c}
2\left( -1\right) ^{n},\text{ if }i=0 \\
\left( 2C_{n}^{i}+C_{n}^{i-1}\right) \left( -1\right) ^{n-i},\text{ if }%
1\leq i\leq n \\
-1,\text{ if }i=n+1.%
\end{array}%
\right.
\end{equation*}%
Since $\deg P=m>1,$ then for any two integers $j$ and $k$ such that $0\leq
j<k\leq n+1,$ we have%
\begin{equation*}
\rho \left( e^{-P_{kc}+P_{jc}}\right) =\deg P-1.
\end{equation*}%
It's clear now that all the conditions of Lemma 2.3 are satisfied. So, by
Lemma 2.3 we obtain $b_{i}\equiv 0$ for all $i=0,...,n+1,$ which is
impossible. Then, $\alpha _{n+1}\left( z\right) \not\equiv 0.$ By Lemma 2.2,
$\left( 3.14\right) $ and $\left( 3.15\right) ,$ we deduce that $\lambda
\left( e^{P}\right) =\deg P>1,$ which is a contradiction.

\quad

\noindent $\left( \text{ii}\right) $ $\deg P=1.$ Suppose now that $P\left(
z\right) =\mu z+\eta $ $\left( \mu \neq 0\right) .$ Assume that $\alpha
_{n+1}\left( z\right) \equiv 0.$ It easy to see that%
\begin{equation*}
\Delta _{c}^{n}\left( 2e^{-P}-e^{-P_{c}}\right) =\left( 2-e^{-\mu c}\right)
\Delta _{c}^{n}\left( e^{-P}\right) .
\end{equation*}%
In the following two subcases, we prove that both of $\left( 2-e^{-\mu
c}\right) $ and $\Delta _{c}^{n}\left( e^{-P}\right) $ are not vanishing.

\noindent \textbf{(A)} Suppose that $2=e^{-\mu c}.$ Then for any integer $i$%
, we have $e^{-i\mu c}=2^{i}$ and $e^{-P_{ic}}=2^{i}e^{-P},$ applying that
on the equation $\left( 3.6\right) $, we\ get%
\begin{equation*}
A_{n}\left( z\right) =\overset{n}{\underset{i=0}{\sum }}C_{n}^{i}\left(
-1\right) ^{n-i}2^{i}e^{-i\mu c}-e^{P}=3^{n}-e^{P},
\end{equation*}%
and
\begin{equation*}
B_{n}\left( z\right) =a\left( z\right) \left( \overset{n}{\underset{i=0}{%
\sum }}C_{n}^{i}\left( -1\right) ^{n-i}\left( 2^{i}-1\right)
e^{-P_{ic}}-1\right)
\end{equation*}%
\begin{equation*}
=a\left( z\right) \left( \overset{n}{\underset{i=0}{\sum }}C_{n}^{i}\left(
-1\right) ^{n-i}\left( 4^{i}-2^{i}\right) e^{-P}-1\right) =a\left( z\right)
\left( \left( 3^{n}-1\right) e^{-P}-1\right) .
\end{equation*}%
Then%
\begin{equation*}
\left( 3^{n}-e^{P}\right) g\left( z\right) +a\left( z\right) \left( \left(
3^{n}-1\right) e^{-P}-1\right) =0,
\end{equation*}%
which is equivalent to%
\begin{equation}
g\left( z\right) =a\left( z\right) \frac{e^{P}-\left( 3^{n}-1\right) }{%
e^{P}\left( 3^{n}-e^{P}\right) }.  \tag{3.17}
\end{equation}%
By the same arguing as before and the equation $\left( 3.7\right) ,$ we
obtain%
\begin{equation*}
g\left( z\right) =a\left( z\right) \frac{e^{P}-\left( 3^{n+1}-1\right) }{%
e^{P}\left( 3^{n+1}-e^{P}\right) },
\end{equation*}%
which contradicts $\left( 3.17\right) $.

\noindent \textbf{(B) }Suppose now that\textbf{\ }$\Delta _{c}^{n}\left(
e^{-P}\right) \equiv 0.$ Thus%
\begin{equation*}
\Delta _{c}^{n}\left( e^{-P}\right) =\overset{n}{\underset{i=0}{\sum }}%
C_{n}^{i}\left( -1\right) ^{n-i}e^{-\mu \left( z+ic\right) -\eta }=e^{-P}%
\overset{n}{\underset{i=0}{\sum }}C_{n}^{i}\left( -1\right) ^{n-i}e^{-\mu ic}
\end{equation*}%
\begin{equation*}
=e^{-P}\left( e^{-\mu c}-1\right) ^{n}.
\end{equation*}%
This together with $\Delta _{c}^{n}e^{-P}\equiv 0$ gives $\left( e^{-\mu
c}-1\right) ^{n}\equiv 0,$ which yields $e^{\mu c}\equiv 1.$ Therefore, for
any $j\in
\mathbb{Z}
$%
\begin{equation}
e^{P\left( z+jc\right) }=e^{\mu z+\mu jc+\eta }=\left( e^{\mu c}\right)
^{j}e^{P\left( z\right) }=e^{P\left( z\right) }.  \tag{3.18}
\end{equation}%
On the other hand, we have from $\left( 3.1\right) $%
\begin{equation}
\Delta _{c}^{n}f\left( z\right) =e^{P\left( z\right) }\left( f\left(
z\right) -a\left( z\right) \right) +a\left( z\right) .  \tag{3.19}
\end{equation}%
By $\left( 3.18\right) $ and $\left( 3.19\right) ,$ we have
\begin{equation}
\Delta _{c}^{n+1}f\left( z\right) =e^{P\left( z\right) }\Delta _{c}f\left(
z\right)  \tag{3.20}
\end{equation}%
Combining $\left( 3.2\right) $ and $\left( 3.20\right) ,$ we obtain
\begin{equation*}
\Delta _{c}f\left( z\right) =\left( f\left( z\right) -a\left( z\right)
\right) +a\left( z\right) e^{-P\left( z\right) }
\end{equation*}%
which means that $\Delta _{c}^{n+1}f\left( z\right) =\Delta _{c}^{n}f\left(
z\right) $ for all $n\geq 1.$ Therefore, \textit{\ }$f\left( z\right) ,$%
\textit{\ }$\Delta _{c}f\left( z\right) $\textit{\ }and $\Delta
_{c}^{2}f\left( z\right) $\ share $a\left( z\right) $\ CM and by Lemma 2.4
we obtain $\Delta _{c}f\left( z\right) =f\left( z\right) ,$ which
contradicts the hypothesis. Then $\Delta _{c}^{n}\left( e^{-P}\right)
\not\equiv 0.$ From the subcases (A) and (B), we can deduce that $\alpha
_{n+1}\left( z\right) \not\equiv 0.$ It is clear that
\begin{equation*}
\max \left\{ \rho \left( a_{i}\right) ,\rho \left( \alpha _{n+1}\right)
,i=0,...,n+1\right\} <\deg P=1.
\end{equation*}%
By using Lemma 2.2, we obtain $\lambda \left( e^{P}\right) =\deg P=1,$ which
is a contradiction, and $P$ must be a constant.

\noindent \textbf{Case 2.} $P\left( z\right) \equiv K$ $,$ $K\in
\mathbb{C}
-\left\{ 2k\pi i,k\in
\mathbb{Z}
\right\} .$ We have from $\left( 3.1\right) $%
\begin{equation*}
\Delta _{c}^{n}f\left( z\right) =e^{K}\left( f\left( z\right) -a\left(
z\right) \right) +a\left( z\right) .
\end{equation*}%
Hence
\begin{equation}
\Delta _{c}^{n+1}f\left( z\right) =e^{K}\Delta _{c}f\left( z\right) .
\tag{3.21}
\end{equation}%
Combining $\left( 3.2\right) $ and $\left( 3.21\right) ,$ we obtain
\begin{equation*}
\Delta _{c}f\left( z\right) =\left( f\left( z\right) -a\left( z\right)
\right) +a\left( z\right) e^{-K}
\end{equation*}%
which means that $\Delta _{c}^{n+1}f\left( z\right) =\Delta _{c}^{n}f\left(
z\right) $ for all $n\geq 1.$ Therefore, \textit{\ }$f\left( z\right) ,$%
\textit{\ }$\Delta _{c}f\left( z\right) $\textit{\ }and $\Delta
_{c}^{2}f\left( z\right) $\ share $a\left( z\right) $\ CM and by Lemma 2.4
we obtain $\Delta _{c}f\left( z\right) =f\left( z\right) ,$ which
contradicts the hypothesis. Then $e^{P}\equiv 1$ and the proof of Theorem
1.1 is completed.

\quad

\noindent \textbf{Proof of the Theorem 1.2. }Setting $g\left( z\right)
=f\left( z\right) +b\left( z\right) -a\left( z\right) .$ Since $\Delta
_{c}^{m}a\left( z\right) \equiv 0$ $\left( 1\leq m\leq n\right) ,$ then we
can remark that
\begin{equation*}
g\left( z\right) -b\left( z\right) =f\left( z\right) -a\left( z\right) ,
\end{equation*}%
\begin{equation*}
\Delta _{c}^{n}g\left( z\right) -b\left( z\right) =\Delta _{c}^{n}f\left(
z\right) -b\left( z\right)
\end{equation*}%
and
\begin{equation*}
\Delta _{c}^{n+1}g\left( z\right) -b\left( z\right) =\Delta _{c}^{n}f\left(
z\right) -b\left( z\right) ,\text{ }n\geq 2.
\end{equation*}%
Since $f\left( z\right) -a\left( z\right) ,$\textit{\ }$\Delta
_{c}^{n}f\left( z\right) -b\left( z\right) $\textit{\ }and $\Delta
_{c}^{n+1}f\left( z\right) -b\left( z\right) $\ share $0$\ CM, then $g\left(
z\right) ,$\textit{\ }$\Delta _{c}^{n}g\left( z\right) $\textit{\ }and $%
\Delta _{c}^{n+1}g\left( z\right) $ share $b\left( z\right) $ CM. By using
Theorem 1.1, we deduce that $\Delta _{c}g\left( z\right) \equiv g\left(
z\right) ,$ which leads to $\Delta _{c}f\left( z\right) \equiv f\left(
z\right) +b\left( z\right) +\Delta _{c}a\left( z\right) -a\left( z\right) $
and the proof of Theorem 1.2 is completed.

\quad

\noindent \textbf{Proof of the Theorem 1.3. }Note that $f\left( z\right) $
is a nonconstant entire function of finite order. Since $f\left( z\right) ,$%
\textit{\ }$\Delta _{c}^{n}f\left( z\right) $\textit{\ }and $\Delta
_{c}^{n+1}f\left( z\right) $\ share $0$\textit{\ }CM, then it is known by
Theorem F that $\Delta _{c}^{n+1}f\left( z\right) =C\Delta _{c}^{n}f\left(
z\right) ,$ where $C$ is a nonzero constant. Then we have%
\begin{equation}
\frac{\Delta _{c}^{n}f\left( z\right) }{f\left( z\right) }=e^{P\left(
z\right) }  \tag{3.22}
\end{equation}%
and
\begin{equation}
\frac{\Delta _{c}^{n+1}f\left( z\right) }{f\left( z\right) }=Ce^{P\left(
z\right) },  \tag{3.23}
\end{equation}%
where $P$ is a polynomial. By $\left( 3.22\right) $ and $\left( 3.23\right) $
we obtain%
\begin{equation}
f_{ic}\left( z\right) =\left( C+1\right) ^{i}e^{P-P_{ic}}f\left( z\right) .
\tag{3.24}
\end{equation}%
Then
\begin{equation}
\Delta _{c}^{n}f\left( z\right) =\left( \overset{n}{\underset{i=0}{\sum }}%
C_{n}^{i}\left( -1\right) ^{n-i}\left( C+1\right) ^{i}e^{P-P_{ic}}\right)
f\left( z\right) =e^{P\left( z\right) }f\left( z\right) .  \tag{3.25}
\end{equation}%
The equality $\left( 3.25\right) $ leads to $\deg P=0.$ Hence $P\left(
z\right) -P_{ic}\left( z\right) \equiv 0$ and $\left( 3.25\right) $ will be%
\begin{equation}
\overset{n}{\underset{i=0}{\sum }}C_{n}^{i}\left( -1\right) ^{n-i}\left(
C+1\right) ^{i}=C^{n}=e^{P\left( z\right) }.  \tag{3.26}
\end{equation}%
By $\left( 3.22\right) ,\left( 3.23\right) $ and $\left( 3.26\right) $ we
deduce
\begin{equation*}
\Delta _{c}^{n}f\left( z\right) =C^{n}f\left( z\right)
\end{equation*}%
and
\begin{equation*}
\Delta _{c}^{n+1}f\left( z\right) =C^{n+1}f\left( z\right) .
\end{equation*}%
Then%
\begin{equation*}
\Delta _{c}^{n+1}f\left( z\right) =\Delta _{c}\left( \Delta _{c}^{n}f\left(
z\right) \right) =\Delta _{c}\left( C^{n}f\left( z\right) \right)
=C^{n}\Delta _{c}f\left( z\right) =C^{n+1}f\left( z\right) ,
\end{equation*}%
which implies $\Delta _{c}f\left( z\right) =Cf\left( z\right) .$ Thus, the
proof of Theorem 1.3 is completed.

\quad

\noindent \textbf{Proof of Corollary 1.2. }By Theorem 1.3 we have $\Delta
_{c}f\left( z\right) =Cf\left( z\right) ,$ where $C$ is a nonzero constant.
Then
\begin{equation}
\Delta _{c}^{m}f\left( z\right) =C\Delta _{c}^{m-1}f\left( z\right)
=C^{m}f\left( z\right) ,\text{ }m\geq 1.  \tag{3.27}
\end{equation}%
On the other hand, for $z_{0}\in
\mathbb{C}
$ we have
\begin{equation}
\Delta _{c}^{m}f\left( z_{0}\right) =f\left( z_{0}\right) .  \tag{3.28}
\end{equation}%
By $\left( 3.27\right) $ and $\left( 3.28\right) $ we deduce that $C^{m}=1.$
Hence $\Delta _{c}^{m}f\left( z\right) =f\left( z\right) $.

\quad

\noindent \textbf{Proof of the Theorem 1.4. }Setting $g\left( z\right)
=f\left( z\right) -a\left( z\right) ,$ we can remark that
\begin{equation*}
g\left( z\right) =f\left( z\right) -a\left( z\right) ,
\end{equation*}%
\begin{equation*}
\Delta _{c}^{n}g\left( z\right) =\Delta _{c}^{n}f\left( z\right) -b\left(
z\right)
\end{equation*}%
and
\begin{equation*}
\Delta _{c}^{n+1}g\left( z\right) =\Delta _{c}^{n}f\left( z\right) -b\left(
z\right) ,\text{ }n\geq 2.
\end{equation*}%
Since $f\left( z\right) -a\left( z\right) ,$\textit{\ }$\Delta
_{c}^{n}f\left( z\right) -b\left( z\right) $\textit{\ }and $\Delta
_{c}^{n+1}f\left( z\right) -b\left( z\right) $\ share $0$\ CM, then $g\left(
z\right) ,$\textit{\ }$\Delta _{c}^{n}g\left( z\right) $\textit{\ }and $%
\Delta _{c}^{n+1}g\left( z\right) $ share $0$ CM. By using Theorem 1.3, we
deduce that $\Delta _{c}g\left( z\right) \equiv Cg\left( z\right) ,$ where $%
C $\ is a nonzero constant\textit{,} which leads to $\Delta _{c}f\left(
z\right) \equiv Cf\left( z\right) +\Delta _{c}a\left( z\right) -a\left(
z\right) $ and the proof of Theorem 1.4 is completed.

\begin{center}
{\Large References}
\end{center}

\noindent $\left[ 1\right] \ $W. Bergweiler, J. K. Langley, \textit{Zeros of
differences of meromorphic functions}, Math. Proc. Cambridge Philos. Soc.
142 (2007), no. 1, 133--147.

\noindent $\left[ 2\right] \ $B. Chen, Z. X. Chen and S. Li,\textit{\
Uniqueness theorems on entire functions and their difference operators or
shifts}, Abstr. Appl. Anal. 2012, Art. ID 906893, 8 pp.

\noindent $\left[ 3\right] \ $B. Chen, and S. Li, \textit{Uniquness problems
on entire functions that share a small function with their difference
operators}, Adv. Difference Equ. 2014, 2014:311, 11 pp.

\noindent $\left[ 4\right] \ $Z. X. Chen,\textit{\ Zeros of entire solutions
to complex linear difference equations, }Acta Math. Sci. Ser. B Engl. Ed. 32
(2012), no. 3, 1141--1148.

\noindent $\left[ 5\right] \ $Y. M. Chiang, S. J. Feng, \textit{On the
Nevanlinna characteristic of }$f\left( z+\eta \right) $ \textit{and
difference equations in the complex plane, }Ramanujan J. 16 (2008), no. 1,
105-129.

\noindent $\left[ 6\right] \ $A. El Farissi, Z. Latreuch and A. Asiri,
\textit{On the uniqueness theory of entire functions and their difference
operators. }Submitted.

\noindent $\left[ 7\right] \ $R. G. Halburd, R. J. Korhonen, \textit{%
Difference analogue of the lemma on the logarithmic derivative with
applications to difference equations, }J. Math. Anal. Appl. 314 (2006)%
\textit{, }no. 2, 477-487.

\noindent $\left[ 8\right] \ $R. G. Halburd, R. J. Korhonen, \textit{%
Nevanlinna theory for the difference operator}, Ann. Acad. Sci. Fenn. Math.
31 (2006), no. 2, 463--478.

\noindent $\left[ 9\right] \ $W. K. Hayman, \textit{Meromorphic functions},
Oxford Mathematical Monographs Clarendon Press, Oxford 1964.

\noindent $\left[ 10\right] \ $G. Jank. E. Mues and L. Volkmann, \textit{%
Meromorphe Funktionen, die mit ihrer ersten und zweiten Ableitung einen
endlichen Wert teilen}, Complex Variables Theory Appl. 6 (1986), no. 1,
51--71.

\noindent $\left[ 11\right] \ $I. Laine, \textit{Nevanlinna theory and
complex differential equations}, de Gruyter Studies in Mathematics, 15.
Walter de Gruyter \& Co., Berlin, 1993.

\noindent $\lbrack 12]$ Z. Latreuch and B. Bela\"{\i}di, \textit{Growth and
oscillation of meromorphic solutions of linear difference equations, }Mat.
Vesnik 66 (2014), no. 2, 213--222.

\noindent $\lbrack 13]$ Z. Latreuch, A. El Farissi and B. Bela\"{\i}di,
\textit{Entire functions sharing small functions with their difference
operators.} Electron. J. Diff. Equ., Vol. 2015 (2015), No. 132, 1-10.

\noindent $\lbrack 14]$ P. Li and C. C. Yang, \textit{Uniqueness theorems on
entire functions and their derivatives}, J. Math. Anal. Appl. 253 (2001),
no. 1, 50--57.

\noindent $\left[ 15\right] $ L.A. Rubel and C. C. Yang, \textit{Values
shared by an entire function and its derivatives}, Lecture Notes in Math.
599(1977), Berlin, Springer - Verlag, 101-103.

\noindent $\left[ 16\right] \ $S. Z. Wu and X. M. Zheng, \textit{Growth of
solutions of some kinds of linear difference equations, }Adv. Difference
Equ. (2015) 2015:142, 11 pp.

\noindent $\left[ 17\right] \ $C. C. Yang, H. X. Yi, \textit{Uniqueness
theory of meromorphic functions}, Mathematics and its Applications, 557.
Kluwer Academic Publishers Group, Dordrecht, 2003.

\noindent $\left[ 18\right] $ L. Z. Yang, \textit{Further results on entire
functions that share one value with their derivatives}, J. Math. Anal. Appl.
212 (1997), 529-536.

\end{document}